\documentclass[10pt]{amsart}

\usepackage{amssymb,amsfonts}

\begin{document}
\theoremstyle{plain}
\newtheorem{thm}{Theorem}
\newtheorem{theorem}[thm]{Theorem}
\newtheorem{lemma}[thm]{Lemma}
\newtheorem{corollary}[thm]{Corollary}
\newtheorem{proposition}[thm]{Proposition}
\newtheorem{conjecture}[thm]{Conjecture}
\theoremstyle{definition}
\newtheorem{construction}[thm]{Construction}
\newtheorem{notations}[thm]{Notations}
\newtheorem{question}[thm]{Question}
\newtheorem{problem}[thm]{Problem}
\newtheorem{remark}[thm]{Remark}
\newtheorem{remarks}[thm]{Remarks}
\newtheorem{definition}[thm]{Definition}
\newtheorem{claim}[thm]{Claim}
\newtheorem{assumption}[thm]{Assumption}
\newtheorem{assumptions}[thm]{Assumptions}
\newtheorem{properties}[thm]{Properties}
\newtheorem{example}[thm]{Example}
\newtheorem{comments}[thm]{Comments}
\newtheorem{blank}[thm]{}
\newtheorem{observation}[thm]{Observation}
\newtheorem{defn-thm}[thm]{Definition-Theorem}


\title[New results of intersection numbers]{New results of intersection numbers\\ on moduli spaces of curves}
        \author{Kefeng Liu}
        \address{Center of Mathematical Sciences, Zhejiang University, Hangzhou, Zhejiang 310027, China;
                Department of Mathematics,University of California at Los Angeles,
                Los Angeles, CA 90095-1555, USA}
        \email{liu@math.ucla.edu, liu@cms.zju.edu.cn}
        \author{Hao Xu}
        \address{Center of Mathematical Sciences, Zhejiang University, Hangzhou, Zhejiang 310027, China}
        \email{haoxu@cms.zju.edu.cn}

        \begin{abstract}
        We present a series of new results we obtained recently about the intersection
        numbers of tautological classes on moduli spaces
        of curves, including a simple formula of the n-point functions for Witten's $\tau$ classes,
        an effective recursion formula to compute
        higher Weil-Petersson volumes, several new recursion formulae of intersection
        numbers and our proof of a conjecture of Itzykson and Zuber
        concerning denominators of intersection numbers. We also present Virasoro and KdV
        properties of generating functions of
        general mixed $\kappa$ and $\psi$ intersections.
        \end{abstract}

    \maketitle

\vskip 10pt {\noindent\bf\large Introduction} \vskip 5pt

Let $\overline{\mathcal M}_{g,n}$ denote the Deligne-Mumford moduli
stack of stable curves of genus $g$ with $n$ marked points. Let
$\psi_i$ be the first Chern class of the line bundle whose fiber
over each pointed stable curve is the cotangent line at the $i$-th
marked point. Let $\lambda_i$ be the $i$-th Chern class of the Hodge
bundle $\mathbb E$, whose fiber over each pointed stable curves is
$H^0(C,\omega_C)$.

We also have the $\kappa$ classes originally defined by Mumford
\cite{Mu}, Morita \cite{Mo} and Miller \cite{Mi}. A more natural
variation was later given by Arbarello-Cornalba \cite{ArCo}. It is
known that the $\kappa$ and $\psi$ classes generate the tautological
cohomology ring of the moduli spaces, and most of the known
cohomology classes are tautological.

The following intersection numbers
$$\langle\tau_{d_1}\cdots\tau_{d_n}\prod_{j\geq1}\kappa_j^{b_j}\rangle_g:=\int_{\overline{\mathcal M}_{g,n}}\psi_1^{d_1}\cdots\psi_n^{d_n}\prod_{j\geq1}\kappa_j^{b_j}.$$
are called the higher Weil-Petersson volumes \cite{KMZ}. These are
important invariants of moduli spaces of curves.

In 1990, Witten \cite{Wi} made the remarkable conjecture that the
generating function of intersection numbers of $\psi$ classes on
moduli spaces are governed by KdV hierarchy. Witten's conjecture
(first proved by Kontsevich \cite{Ko}) is among the deepest known
properties of moduli spaces of curves and motivated a surge of
subsequent developments.

The intersection theory of tautological classes on the moduli space
of curves is a very important subject and has close connections to
string theory, quantum gravity and many branches of mathematics.

\vskip 10pt {\noindent\bf\large The n-point functions for
intersection numbers}\vskip 5pt

\begin{definition}
We call the following generating function
$$F(x_1,\dots,x_n)=\sum_{g=0}^{\infty}\sum_{\sum d_j=3g-3+n}\langle\tau_{d_1}\cdots\tau_{d_n}\rangle_g\prod_{j=1}^n x_j^{d_j}$$
the $n$-point function.
\end{definition}

Consider the following ``normalized'' $n$-point function
$$G(x_1,\dots,x_n)=\exp\left(\frac{-\sum_{j=1}^n
x_j^3}{24}\right)\cdot F(x_1,\dots,x_n).$$ Starting from $1$-point
function $G(x)=\frac{1}{x^2}$, we can obtain any $n$-point function
recursively by the following theorem.

\begin{theorem} \cite{LX3} For $n\geq2$,
$$G(x_1,\dots,x_n)=\sum_{r,s\geq0}\frac{(2r+n-3)!!}{4^s(2r+2s+n-1)!!}P_r(x_1,\dots,x_n)\Delta(x_1,\dots,x_n)^s,$$
where $P_r$ and $\Delta$ are homogeneous symmetric polynomials
defined by
\begin{align*}
\Delta(x_1,\dots,x_n)&=\frac{(\sum_{j=1}^nx_j)^3-\sum_{j=1}^nx_j^3}{3},\nonumber\\
P_r(x_1,\dots,x_n)&=\left(\frac{1}{2\sum_{j=1}^nx_j}\sum_{\underline{n}=I\coprod
J}(\sum_{i\in I}x_i)^2(\sum_{i\in J}x_i)^2G(x_I)
G(x_J)\right)_{3r+n-3}\nonumber\\
&=\frac{1}{2\sum_{j=1}^nx_j}\sum_{\underline{n}=I\coprod
J}(\sum_{i\in I}x_i)^2(\sum_{i\in J}x_i)^2\sum_{r'=0}^r
G_{r'}(x_I)G_{r-r'}(x_J),
\end{align*}
where $I,J\ne\emptyset$, $\underline{n}=\{1,2,\ldots,n\}$ and
$G_g(x_I)$ denotes the degree $3g+|I|-3$ homogeneous component of
the normalized $|I|$-point function $G(x_{k_1},\dots,x_{k_{|I|}})$,
where $k_j\in I$.
\end{theorem}

Thus we have an elementary and more efficient algorithm to calculate
all intersection numbers of $\psi$ classes other than the celebrated
Witten-Kontsevich theorem.

Since $P_0(x,y)=\frac{1}{x+y}$, $P_r(x,y)=0$ for $r>0$ and
$$P_r(x,y,z)=\frac{r!}{2^r(2r+1)!}\frac{(xy)^r(x+y)^{r+1}+(yz)^r(y+z)^{r+1}+(zx)^r(z+x)^{r+1}}{x+y+z},$$
we recover Dijkgraaf's $2$-point function and Zagier's $3$-point
function obtained more than ten years ago.

There is another slightly different formula of the $n$-point
functions. When $n=3$, this has also been obtained by Zagier.
\begin{theorem} \cite{LX3} For $n\geq2$,
$$
F(x_1,\dots,x_n)=\exp\frac{(\sum_{j=1}^n
x_j)^3}{24}\sum_{r,s\geq0}\frac{(-1)^s
P_r(x_1,\dots,x_n)\Delta(x_1,\dots,x_n)^s}{8^s(2r+2s+n-1)s!}
$$
where $P_r$ and $\Delta$ are the same polynomials as defined in
Theorem 2.
\end{theorem}

Okounkov \cite{Ok} obtained an analytic expression of the $n$-point
functions using $n$-dimensional error-function-type integrals.
Br\'ezin and Hikami \cite{BH} use correlation functions of GUE
ensemble to find explicit formulae of $n$-point functions.

\vskip 10pt {\noindent\bf\large Recursion formulae of higher
Weil-Petersson volumes} \vskip 5pt

We have discovered a general recursion formula of higher
Weil-Petersson volumes \cite{LX5}, which is a vast generalization of
the Mirzakhani recursion formula \cite{LX4}.

First we fix notations as in \cite{KMZ}.

Consider the semigroup $N^\infty$ of sequences ${\bold
m}=(m_1,m_2,\dots)$ where $m_i$ are nonnegative integers and $m_i=0$
for sufficiently large $i$.

Let $\bold m, \bold t, \bold{a_1,\dots,a_n} \in N^\infty$, $\bold
m=\sum_{i=1}^n \bold{a_i}$, and $\bold s:=(s_1,s_2,\dots)$ be a
family of independent formal variables.
$$|\bold m|:=\sum_{i\geq 1}i m_i,\quad ||\bold m||:=\sum_{i\geq1}m_i,\quad \bold s^{\bold m}:=\prod_{i\geq 1}s_i^{m_i},\quad \bold m!:=\prod_{i\geq1}m_i!,$$
$$\binom{\bold m}{\bold{t}}:=\prod_{i\geq1}\binom{m_i}{t_i},\quad \binom{\bold m}{\bold{a_1,\dots,a_n}}:=\prod_{i\geq1}\binom{ m_i}{a_1(i),\dots,a_n(i)}.$$

Let $\bold b\in N^\infty$, we denote a formal monomial of $\kappa$
classes by
$$\kappa(\bold b):=\prod_{i\geq1}\kappa_i^{b_i}.$$

\begin{theorem} \cite{LX5} Let $\bold b\in N^\infty$ and $d_j\geq 0$.
\begin{multline*}
(2d_1+1)!!\langle\kappa(\bold b)\prod_{j=1}^n\tau_{d_j}\rangle_g\\
=\sum_{j=2}^n\sum_{\bold L+\bold{L'}=\bold b}\alpha_{\bold
L}\binom{\bold b}{\bold L}\frac{(2(|\bold
L|+d_1+d_j)-1)!!}{(2d_j-1)!!}\langle\kappa(\bold{L'})\tau_{|\bold
L|+d_1+d_j-1}\prod_{i\neq
1,j}\tau_{d_i}\rangle_g\\
+\frac{1}{2}\sum_{\bold L+\bold{L'}=\bold b}\sum_{r+s=|\bold
L|+d_1-2}\alpha_{\bold
L}\binom{\bold b}{\bold L}(2r+1)!!(2s+1)!!\langle\kappa(\bold{L'})\tau_r\tau_s\prod_{i\neq1}\tau_{d_i}\rangle_{g-1}\\
+\frac{1}{2}\sum_{\substack{\bold L+\bold{e}+\bold{f}=\bold
b\\I\coprod J=\{2,\dots,n\}}}\sum_{r+s=|\bold L|+d_1-2}\alpha_{\bold
L}\binom{\bold b}{\bold
L,\bold{e},\bold{f}}(2r+1)!!(2s+1)!!\\
\times \langle\kappa(\bold{e})\tau_r\prod_{i\in
I}\tau_{d_i}\rangle_{g'}\langle\kappa(\bold{f})\tau_s\prod_{i\in
J}\tau_{d_i}\rangle_{g-g'}.
\end{multline*}
These tautological constants $\alpha_{\bold L}$ can be determined
recursively from the following formula
$$\sum_{\bold L+\bold{L'}=\bold b}\frac{(-1)^{||\bold L||}\alpha_{\bold L}}{\bold L!\bold{L'}!(2|\bold{L'}|+1)!!}=0,\qquad \bold b\neq0,$$
namely
$$\alpha_{\bold b}=\bold b!\sum_{\substack{\bold L+\bold{L'}=\bold b\\ \bold{L'}\neq\bold 0}}\frac{(-1)^{||\bold L'||-1}\alpha_{\bold L}}{\bold L!\bold{L'}!(2|\bold{L'}|+1)!!},\qquad\bold b\neq0,$$
with the initial value $\alpha_{\bold 0}=1$.
\end{theorem}

The proof of the above theorem is to use Witten-Kontsevich theorem,
a combinatorial formula in \cite{KMZ} expressing $\kappa$ classes by
$\psi$ classes and the following elementary but crucial lemma
\cite{LX5}.
\begin{lemma}
Let $F(\bold L,n)$ and $G(\bold L,n)$ be two functions defined on
$N^\infty\times\mathbb N$, where $\mathbb N=\{0,1,2,\dots\}$ is the
set of nonnegative integers. Let $\alpha_{\bold L}$ and
$\beta_{\bold L}$ be real numbers depending only on $\bold L\in
N^\infty$ that satisfy $\alpha_{\bold 0}\beta_{\bold 0}=1$ and
$$\sum_{\bold L+\bold{L'}=\bold b}\alpha_{\bold L}\beta_{\bold{L'}}=0,\qquad\bold b\neq0.$$
Then the following two identities are equivalent.

\begin{align*}
G(\bold b,n)=\sum_{\bold L+\bold{L'}=\bold b} \alpha_{\bold L} F(\bold{L'},n+|\bold L|),\quad\forall\ (\bold b,n)\in N^\infty\times\mathbb N\\
F(\bold b,n)=\sum_{\bold L+\bold{L'}=\bold b} \beta_{\bold L} G(\bold{L'},n+|\bold L|),\quad\forall\ (\bold b,n)\in N^\infty\times\mathbb N\\
\end{align*}
\end{lemma}

When $\bold b=(l,0,0,\dots)$, Theorem 4 recovers Mirzakhani's
recursion formula of Weil-Petersson volumes for moduli spaces of
bordered Riemann surfaces \cite{Mir,Mir2,MS,Saf}.

Theorem 4 also provides an effective algorithm to compute higher
Weil-Petersson volumes recursively.

In fact we can use the main formula in \cite{KMZ} to generalize
almost all pure $\psi$ intersections to identities of higher
Weil-Petersson volumes which share similar structures as Theorem 4.
For example, the identities in the following theorem are
generalizations of the string and dilaton equations.
\begin{theorem} \cite{LX5}
For $\bold b\in N^\infty$ and $d_j\geq0$,
$$\sum_{\bold L+\bold L'=\bold
b}(-1)^{||\bold L||}\binom{\bold b}{\bold L}\langle\tau_{|\bold
L|}\prod_{j=1}^n\tau_{d_j}\kappa(\bold
L')\rangle_g=\sum_{j=1}^n\langle\tau_{d_j-1}\prod_{i\neq
j}\tau_{d_i}\kappa(\bold b)\rangle_g,$$ and
$$\sum_{\bold L+\bold L'=\bold b}(-1)^{||\bold
L||}\binom{\bold b}{\bold L}\langle\tau_{|\bold
L|+1}\prod_{j=1}^n\tau_{d_j}\kappa(\bold
L')\rangle_g=(2g-2+n)\langle\prod_{j=1}^n\tau_{d_j}\kappa(\bold
b)\rangle_g.$$
\end{theorem}
Note that Theorem 6 generalizes the results in \cite{DN}.

\vskip 10pt {\noindent\bf\large New identities of intersection
numbers} \vskip 5pt

The next two theorems follow from a detailed study of coefficients
of the $n$-point functions in Theorem 2.
\begin{theorem} \cite{LX3} We have
\begin{enumerate}
\item Let $k>2g$, $d_j\geq0$ and $\sum_{j=1}^n
d_j=3g+n-k$.
$$\sum_{\underline{n}=I\coprod
J}\sum_{j=0}^{k}(-1)^j\langle\tau_{j}\tau_0^2\prod_{i\in
I}\tau_{d_i}\rangle_{g'}\langle\tau_{k-j}\tau_0^2\prod_{i\in
J}\tau_{d_i}\rangle_{g-g'}=0.$$

\item Let $d_j\geq1$ and $\sum_{j=1}^n
d_j=g+n$.
$$\sum_{\underline{n}=I\coprod
J}\sum_{j=0}^{2g}(-1)^j\langle\tau_{j}\tau_0^2\prod_{i\in
I}\tau_{d_i}\rangle_{g'}\langle\tau_{2g-j}\tau_0^2\prod_{i\in
J}\tau_{d_i}\rangle_{g-g'}=\frac{(2g+n+1)!}{4^g(2g+1)!\prod_{j=1}^n(2d_j-1)!!}.$$
\end{enumerate}
\end{theorem}

\begin{theorem} \cite{LX2,LX3} We have
\begin{enumerate}

\item Let $k>2g$, $d_j\geq0$ and $\sum_{j=1}^{n}d_j=3g+n-k-1$.
$$
\sum_{j=0}^k(-1)^j\langle\tau_{k-j}\tau_j\prod_{i=1}^n\tau_{d_i}\rangle_g=0.
$$

\item Let $d_j\geq1$ and $\sum_{j=1}^{n}(d_j-1)=g-1$.
$$\sum_{j=0}^{2g}(-1)^j\langle\tau_{2g-j}\tau_j\prod_{i=1}^n\tau_{d_i}\rangle_g=\frac{(2g+n-1)!}{4^g(2g+1)!\prod_{j=1}^n(2d_j-1)!!}.$$
\end{enumerate}
\end{theorem}

In fact, it's easy to see that Theorems 7 and 8 imply each other
through the following proposition.

\begin{proposition} \cite{LX3}
Let $d_j\geq0$ and $\sum_{j=1}^n d_j=g+n$.
\begin{multline*}
\sum_{\underline{n}=I\coprod
J}\sum_{j=0}^{2g}(-1)^j\left(\langle\tau_{j}\tau_0^2\prod_{i\in
I}\tau_{d_i}\rangle\langle\tau_{2g-j}\tau_0^2\prod_{i\in
J}\tau_{d_i}\rangle+\langle\tau_{j}\tau_{2g-j}\tau_0^2\prod_{i\in
I}\tau_{d_i}\rangle\langle\tau_0^2\prod_{i\in
J}\tau_{d_i}\rangle\right)\\
=(2g+n+1)\sum_{j=0}^{2g}(-1)^j\langle\tau_0\tau_j\tau_{2g-j}\prod_{i=1}^n\tau_{d_i}\rangle_g
\end{multline*}
\end{proposition}

Since ${\rm ch}_k(\mathbb E)=0$ for $k>2g$,
$\lambda_{g}\lambda_{g-1}=(-1)^{g-1}(2g-1)!{\rm ch}_{2g-1}(\mathbb
E)$, by Mumford's formula \cite{Mu} of the Chern character of Hodge
bundles, it's not difficult to see that Theorem 8 implies the
following theorem.
\begin{theorem} \cite{LX2,LX3}
Let $k$ be an even number and $k\geq 2g$, $d_j\geq0$,
$\sum_{j=1}^{n}d_j=3g+n-k-2$.
\begin{multline}
\langle\prod_{j=1}^n\tau_{d_j}\tau_k\rangle_g=\sum_{j=1}^{n}
\langle\tau_{d_j+k-1}\prod_{i\neq j}\tau_{d_i}\rangle_g\\
-\frac{1}{2}\sum_{\underline{n}=I\coprod
J}\sum_{j=0}^{k-2}(-1)^j\langle\tau_{j}\prod_{i\in
I}\tau_{d_i}\rangle_{g'}\langle\tau_{k-2-j}\prod_{i\in
J}\tau_{d_i}\rangle_{g-g'}.
\end{multline}
\end{theorem}
Note that when $k=2g$, the above theorem is equivalent to the
following Hodge integral identity \cite{GP} (also known as Faber's
intersection number conjecture \cite{Fab})
$$\int_{\overline{\mathcal
M}_{g,n}}\psi_1^{d_1}\dots\psi_n^{d_n}\lambda_g\lambda_{g-1}=\frac{(2g-3+n)!|B_{2g}|}{2^{2g-1}(2g)!\prod_{j=1}^{n}(2d_j-1)!!}.
$$
where $\sum_{j=1}^n (d_j-1)=g-2$ and $d_j\geq1$.

The above $\lambda_g\lambda_{g-1}$ integral follows from degree $0$
Virasoro constraints for $\mathbb P^2$ announced by Givental
\cite{Giv}. However it is very desirable to have a direct proof of
identity (1) when $k=2g$, possibly using our explicit formulae of
the $n$-point functions (see also \cite{GJV}).

As pointed out in the last section, we can generalize all of the
above new recursion formulae of $\psi$ classes to identities of
higher Weil-Petersson volumes. For example, we may generalize
Proposition 9 and Theorem 10 to the following
\begin{proposition} \cite{LX5}
Let $\bold b\in N^\infty$, $d_j\geq0$.
\begin{multline*}
\sum_{j=0}^{2g}(-1)^j\langle\tau_0\tau_1\tau_j\tau_{2g-j}\prod_{j=1}^n\tau_{d_j}\kappa(\bold b)\rangle_g\\
=\sum_{\substack{\bold L+\bold{L'}=\bold b\\
\underline{n}=I\coprod J}}\sum_{j=0}^{2g}(-1)^j\binom{\bold b}{\bold
L}\left(\langle\tau_{j}\tau_0^2\prod_{i\in I}\tau_{d_i}\kappa(\bold
L)\rangle\langle\tau_{2g-j}\tau_0^2\prod_{i\in
J}\tau_{d_i}\kappa(\bold{L'})\rangle\right.\\
\left.+\langle\tau_{j}\tau_{2g-j}\tau_0^2\prod_{i\in
I}\tau_{d_i}\kappa(\bold L)\rangle\langle\tau_0^2\prod_{i\in
J}\tau_{d_i}\kappa(\bold{L'})\rangle\right)\\
\end{multline*}
\end{proposition}

\begin{theorem} \cite{LX5}
Let $\bold b\in N^\infty$, $M\geq 2g$ be an even number and
$d_j\geq0$.
\begin{multline*}
\sum_{\bold L+\bold L'=\bold b}(-1)^{||\bold L||}\binom{\bold
b}{\bold L}\langle\tau_{|\bold
L|+M}\prod_{j=1}^n\tau_{d_j}\kappa(\bold{
L'})\rangle_g=\sum_{j=1}^n\langle\tau_{d_j+M-1}\prod_{i\neq
j}\tau_{d_i}\kappa(\bold b)\rangle_g\\
-\frac{1}{2}\sum_{\substack{\bold L+\bold{L'}=\bold
b\\\underline{n}=I\coprod J}}\sum_{j=0}^{M-2}(-1)^j\binom{\bold
b}{\bold L}\langle\tau_j\prod_{i\in I}\tau_{d_i}\kappa(\bold
L)\rangle_{g'}\langle\tau_{M-2-j}\prod_{i\in
J}\tau_{d_i}\kappa(\bold{L'})\rangle_{g-g'}.
\end{multline*}
\end{theorem}

We also found the following conjectural identity experimentally,
which is amazing if compared with Theorems 8 and 10.
\begin{conjecture} \cite{LX2}
Let $g\geq2$, $d_j\geq1$, $\sum_{j=1}^{n}(d_j-1)=g$.
\begin{multline*}
\frac{(2g-3+n)!}{2^{2g+1}(2g-3)!\prod_{j=1}^n(2d_j-1)!!}=\langle\prod_{j=1}^n\tau_{d_j}\tau_{2g-2}\rangle_g-\sum_{j=1}^{n}
\langle\tau_{d_j+2g-3}\prod_{i\neq j}\tau_{d_i}\rangle_g\\
+\frac{1}{2}\sum_{\underline{n}=I\coprod
J}\sum_{j=0}^{2g-4}(-1)^j\langle\tau_{j}\prod_{i\in
I}\tau_{d_i}\rangle_{g'}\langle\tau_{2g-4-j}\prod_{i\in
J}\tau_{d_i}\rangle_{g-g'}.
\end{multline*}
\end{conjecture}

Since $(2g-3)!{\rm ch}_{2g-3}(\mathbb
E)=(-1)^{g-1}(3\lambda_{g-3}\lambda_g-\lambda_{g-1}\lambda_{g-2})$,
it's easy to see that the above identity is equivalent to the
following identity of Hodge integrals.
\begin{conjecture}
Let $g\geq2$, $d_j\geq1$, $\sum_{j=1}^{n}(d_j-1)=g$.
\begin{multline*}
\frac{2g-2}{|B_{2g-2}|}\left(\langle\prod_{j=1}^n\tau_{d_j}\mid\lambda_{g-1}\lambda_{g-2}\rangle_g-3\langle\prod_{j=1}^n\tau_{d_j}\mid\lambda_{g-3}\lambda_{g}\rangle_g\right)\nonumber\\
=\frac{1}{2}\sum_{j=0}^{2g-4}(-1)^j\langle\tau_{2g-4-j}\tau_j\prod_{i=1}^n\tau_{d_i}\rangle_{g-1}+\frac{(2g-3+n)!}{2^{2g+1}(2g-3)!\prod_{j=1}^n(2d_j-1)!!}.
\end{multline*}
\end{conjecture}

\vskip 10pt {\noindent\bf\large Virasoro constraints and KdV
hierarchy} \vskip 5pt

From Theorem 4, we found new Virasoro constraints and KdV hierarchy
for generating functions of higher Weil-Petersson volumes which
vastly generalize the Witten conjecture and the results of Mulase
and Safnuk \cite{MS}.

Let $\bold s:=(s_1,s_2,\dots)$ and $\bold t:=(t_0,t_1,t_2,\dots)$,
we introduce the following generating function
$$G(\bold s,\bold t):=\sum_{g}\sum_{\bold m,\bold n}\langle\kappa_1^{m_1}\kappa_2^{m_2}\cdots\tau_0^{n_0}\tau_1^{n_1}\cdots\rangle_g\frac{\bold s^{\bold m}}{\bold m!}\prod_{i=0}^\infty\frac{t_i^{n_i}}{n_i!},$$
where $\bold s^{\bold m}=\prod_{i\geq1}s_i^{m_i}$.

We introduce the following family of differential operators for
$k\geq -1$,
 \begin{multline*}
    V_k = -\frac{1}{2} \sum_{\bold L} (2(|\bold L|+k)+3)!! \gamma_{\bold L} \bold
    s^{\bold L}
        \frac{\partial }{\partial t_{|\bold L|+k+1} }
        + \frac{1}{2} \sum_{j=0}^{\infty} \frac{(2(j+k)+1)!! }{(2j-1)!! } t_j
        \frac{\partial }{\partial t_{j+k} } \\
    + \frac{1}{4} \sum_{d_1 + d_2 = k-1}
        (2d_1 + 1)!! (2d_2 + 1)!! \frac{\partial^2 }{\partial t_{d_1} \partial t_{d_2}}
        + \frac{\delta_{k,-1}t_0^2}{4} + \frac{\delta_{k,0} }{48},
 \end{multline*}
 where $\gamma_{\bold L}$ are defined by
 $$\gamma_{\bold L}=\frac{(-1)^{||\bold L||}}{\bold L!(2|\bold
 L|+1)!!}.$$

\begin{theorem} \cite{LX5, MS} We have $V_k\exp(G)=0$ for $k\geq-1$ and
the operators $V_k$ satisfy the Virasoro relations
$$[V_n,V_m]=(n-m)V_{n+m}.$$
\end{theorem}

The Witten-Kontsevich theorem states that the generating function
for $\psi$ class intersections
$$F(t_0, t_1, \ldots)= \sum_{g} \sum_{\bold n} \langle\prod_{i=0}^\infty \tau_{i}^{n_i}\rangle_{g} \prod_{i=0}^\infty \frac{t_i^{n_i} }{n_i!
        }$$
is a $\tau$-function for the KdV hierarchy.

Since Virasoro constraints uniquely determine the generating
functions $G(\bold s, t_0,t_1,\dots)$ and $F(t_0,t_1,\dots)$, we
have the following theorem.
\begin{theorem} \cite{LX5, MS}
$$G(\bold s,t_0,t_1,\dots)=F(t_0,t_1,t_2+p_2,t_3+p_3,\dots),$$
where $p_k$ are polynomials in $\bold s$ given by
$$p_k=\sum_{|\bold L|=k-1}\frac{(-1)^{||\bold L||-1}}{\bold L!}\bold
s^{\bold L}.$$ In particular, for any fixed values of $\bold s$,
$G(\bold s,\bold t)$ is a $\tau$-function for the KdV hierarchy.
\end{theorem}

Theorem 16 also generalized results in \cite{MZ}.

\vskip 10pt {\noindent\bf\large Denominators of intersection
numbers} \vskip 5pt

Let {\it denom}$(r)$ denotes the denominator of a rational number
$r$ in reduced form (coprime numerator and denominator, positive
denominator). We define
$$D_{g,n}=lcm\left\{denom\left(\langle\prod_{j=1}^n\tau_{d_j}\rangle_g\right)\Big{|}\
\sum_{j=1}^{n}d_j=3g-3+n\right\}$$ and for $g\geq 2$,
$$\mathcal D_g=lcm\left\{denom\left(\int_{\overline{\mathcal M}_{g}}\kappa(\bold b)\right)\Big{|}\ |\bold b|=3g-3\right\}$$
where {\it lcm} denotes {\it least common multiple}.

Since denominators of intersection numbers on $\overline{\mathcal
M}_{g,n}$ all come from orbifold quotient singularities, the
divisibility properties of $D_{g,n}$ and $\mathcal D_g$ should
reflect overall behavior of singularities.

We have the following properties of $D_{g,n}$ and $\mathcal D_g$.
\begin{proposition} \cite{LX1}
We have $D_{g,n}\mid D_{g,n+1}$, $D_{g,n} \mid \mathcal D_g$ and
$\mathcal D_g=D_{g,3g-3}$.
\end{proposition}

\begin{theorem} \cite{LX1}
For $1<g'\le g$, the order of any automorphism group of a Riemann
surface of genus $g'$ divides $D_{g,3}$.
\end{theorem}

The following corollary of Theorem 18 is a conjecture raised by
Itzykson and Zuber \cite{IZ} in 1992.
\begin{corollary}
For $1<g'\le g$, the order of any automorphism group of an algebraic
curve of genus $g'$ divides $\mathcal D_g$.
\end{corollary}

The proof of Theorem 18 needs the following two lemmas (see
\cite{LX1}).
\begin{lemma}
If $p\leq g+1$ is a prime number, then ${\rm ord}(p,D_{g,3}) \geq
2$.
\end{lemma}
\begin{lemma} \cite{Ha}
Let $X$ be a Riemann Surface of genus $g\geq2$, then for any prime
number $p$,
$${\rm ord}(p,|Aut(X)|)\leq\lfloor\log_p\frac{2pg}{p-1}\rfloor+{\rm ord}(p,2(g-1)).$$
\end{lemma}

We have also obtained conjectural exact values of $\mathcal D_g$ for
all $g$ in \cite{LX2}.

$$ \ \ \ \ $$

\end{document}